\documentclass[a4paper,11pt]{article}
\usepackage[english]{babel}
\usepackage{amsfonts,makeidx,amssymb,amsmath,amsthm,latexsym}
\usepackage[latin1]{inputenc}
\usepackage[dvips]{graphics}

\newtheorem{ejm}{Example}

\newtheorem{teor}{Theorem}

\newcommand{\C}{\mathbb{C}}

\newcommand{\tx}{\textrm}
\title{On Polars of Plane Branches}
\author{A. Hefez (1), M. E. Hernandes (2), M. F. Hernández Iglesias (2) \\ \\
({\small (1) Universidade Federal Fluminense, (2) Universidade Estadual de Maring\'a })}
\date{ \ }
\begin{document}
\maketitle
It is well known that the equisingularity class of the general polar of a plane branch is not the same for all branches in a given equisingularity class, but it is the same for 
sufficiently general ones (cf. \cite{C1}) and depends upon the analytic type of the branch. The aim of this paper is to go beyond generality and show how one could describe the equisingularity classes of (general) polars of all branches in a given equisingularity class, making use of the analytic classification of branches as described in \cite{HH3}. We will show how this works in some particular equisingularity classes for which one has the complete explicit analytic classification, and in particular for all branches of multiplicity less or equal than four, based on the classification given in \cite{HH2}. 

\section{Introduction}

Let $f\in \C \{x,y\}$ be a convergent power series in two variables over the field of complex numbers with no multiple factors and such that $f(0,0)=0$. Notice that because of finite determinacy the same results hold without changes if the power series are only formal. We consider the germ of curve $(f) \colon f=0$ at the origin of $\C^2$, which determines $f$ up to multiplication by a unit. An irreducible germ of curve will be called a {\em branch}. We will say that the germs of curves $(f)$ and $(g)$ are {\em analytically equivalent} if there exists a germ of analytic diffeomorphism at the origin of $\C^2$ that transforms one germ into the other. In terms of equations, this translates into the fact that there exist a unit $u$ and an automorphism $\phi$ of $\C \{x,y\}$ such that $g=u\phi(f)$. In this case, we also say that the functions $f$ and $g$ are {\em contact equivalent}. The {\em equisingularity class} of a curve is its equivalence class under transformations by germs of homeomorphisms at the origin of $\C^2$.

The polar curve of $f$ in the direction $(a\colon b)\in {\mathbb P}_{\C}^1$ is the germ of curve defined by the equation $af_x+bf_y=0$. It is known (cf. \cite{C2} Theorem 7.2.10) that, except for a finite set of directions, the polar is reduced and its equisingularity class is constant, although its analytic type depends essentially upon the direction $(a:b)$, as we will see in an example at the end of the paper. Also, the equisingularity class of the general polar of $f$ is constant in the contact class of $f$ (cf. \cite{C2} Corollary 8.5.8), but it is not constant in the equisingularity class of $f$, as one can easily check by considering for example the curves $y^3-x^{11}$ and $y^3-x^{11}+x^8y$ (cf. \cite{P}). So, the topological type of the polar of a given curve is not determined only by the topological type of the curve, but it is determined by its analytical type. In the next section we will see to what extent the analytic type of the curve will influence the topology of its polar. 

We refer to \cite{Z} for the definitions and basic results we will use in the sequel. It is a classical result that the equisingularity class of a reduced curve given by $f=f_1\cdots f_r$, where th $f_i$ are irreducible, is determined by the semigroups of the $f_i$'s and their mutual intersection numbers ${\rm I}(f_i,f_j)$, for $i\neq j$. A semigroup of values $\Gamma$ of a branch will be given by its minimal set of generators $ \Gamma=\langle v_0,v_1,\ldots,v_g \rangle$ and the integer $g$ will be called the {\em genus} of the branch. Such a semigroup has a conductor $c$ and the equisingularity class it determines may be parametrized by a constructible set $\mathcal E$ in $\C^{c-v_1-1}$, whose points are the coefficients of the Newton-Puiseux parametrization 
\[
x(t)=t^{v_0}, \  \ y(t)=t^{v_1}+ \sum_{i=v_1+1}^{c-1} c_it^i,
\]
in the sense that any element in the equisingularity class is analytically equivalent to one with a Newton-Puiseux parametrization as above.

Given an equisingularity class of irreducible curves, it was proved in \cite{HH1} and \cite{HH3} that the parameter space $\mathcal E$ may be decomposed into a finite union of disjoint constructible sets $\mathcal E={\mathcal E}_1 \cup \cdots \cup {\mathcal E}_r$, where on each ${\mathcal E}_{\ell}$, the set $\Lambda_\ell$ of values of K\"ahler differentials on the corresponding curve, which is an analytic invariant of the curve, is fixed. 

Since $\Gamma$ has a conductor and $\Gamma \setminus \{0\} \subset \Lambda_\ell$, the set $\Lambda_\ell$ is determined by the finite set $\Lambda_\ell\setminus\Gamma$. If this last set is not empty, the natural number $\lambda$ associated to a curve represented by a point of ${\mathcal E}_\ell$, defined as
\[
\lambda =\min\left( \Lambda_\ell\setminus\Gamma\right)-v_0,
\]
is an analytic invariant known as the {\em Zariski invariant} of the curve.

We will now recall a result that shows that the elements of ${\mathcal E}_\ell$ admit a normal form.

\noindent {\bf Normal Forms Theorem} (cf. \cite{HH3}) {\it If $\mathcal C$ is a curve corresponding to a point in ${\mathcal E}_\ell$, then either $\mathcal C$ is analytically equivalent to a curve with para\-me\-tri\-za\-tion $(t^{v_0},t^{v_1})$, when $\Lambda_\ell \setminus \Gamma =\emptyset$, or to a curve with a parametrization of the form  
\[
x=t^{v_0}, \quad y=t^{v_1}+t^\lambda +\sum_{i}c_i t^i,
\]
where the summation is over all indices $i$ greater than $\lambda$ and do not belong to the set $\Lambda_\ell-v_0$. Moreover, two curves $\mathcal C$, with a parametrization as above, and  ${\mathcal C}'$ with a similar parametrization but with coefficients $(c_i')$ instead of $(c_i)$, are analytically equivalent if and only if there exists a complex number $\zeta$ such that $\zeta^{\lambda-v_1}=1$ and for all $i$, one has $c_i=\zeta^{i-v_1}c_i'$.}\medskip

At this point it is natural to ask if the equisingularity class of the general polar is constant on each stratum ${\mathcal E}_\ell$. 

We will see in the next section that the answer may be negative, but positive for a general member of each irreducible component of the family. This was shown in \cite{C1} in the particular case of the whole space $\mathcal E$.\smallskip

For the convenience of the reader, we will state a well known result about Newton non-degenerate plane curve singularities that will be needed in our analysis. 

We say that a reduced curve $(f)$, where $f=\sum_{ij}a_{ij}x^iy^j$ is {\em Newton non-degenerate} if $x$ and $y$ do not divide $f$ and for any side $L$ of its Newton polygon, the polynomial $f_L=\sum_{(i,j)\in L}a_{ij}x^iy^j$ has no critical points outside the curve $xy=0$ (cf. \cite{K}). Since $f_L$ is a quasi-homogeneous polynomial, we may rephrase the Newton non-degeneration as follows:\smallskip 

Let $P_k=(i_k,j_k)$ and $P_{k+1}=(i_{k+1},j_{k+1})$ be the extremal points of $L$ and define 
$$ p_L(z)=z^{-j_{k+1}}f_L(1,z).
$$
Then one has that $f_L$ has no critical points outside the curve $xy=0$ if and only if $p_L$ has no multiple roots.\smallskip

The following result (cf. \cite{O} or \cite{BLP}) will describe the equisingularity class of a Newton non-degenerate curve $(g)$ such that $(x)$ does not belong to its tangent cone:\medskip

\noindent {\it There is a decomposition $(g_i)$, $i=1,\ldots,s$, of $(g)$, such that the Newton polygon of $(g_i)$ is the segment $[(0,n_i); (m_i,0)]$, with $n_i=I(g_i,x)$ and $m_i=I(g_i,y)$ and  $1 \leq d_1 <  d_2 <  \cdots < d_s \leq \infty$, where $d_i=\frac{m_i}{n_i}$, and $d_s= \infty$ if $(g_s)=(y)$. When  $(g)$ is Newton non-degenerate, then $(g_i)$ has $r_i=\tx{gcd}(n_i,m_i)$ branches given by Newton-Puiseux parametrizations 
$$(g_i)_j \colon y_{i,j}^\frac{n_i}{r_i}=a_{i,j}x^{\frac{m_i}{r_i}}+ \cdots , \ \, 1\leq j \leq r_i, \ \ {\rm with} \ \ a_{i,j}\neq a_{i,j'}, \ \ {\rm for} \ \ j \neq j',$$
and}
$$
{\rm I}((g_i)_j,(g_{i'})_{j'})=\inf\{d_i,d_{i'}\}{\rm I}((g_{i'})_{j'},x){\rm I}((g_{i})_{j},x).
$$

This paper contains some results from the PhD thesis of the third author under the supervision of the other two (cf. \cite{HI}).

\section{Polars and Normal Forms}

We will see in the following example how the Normal Forms Theorem may be used to describe the equisingularity classes of the general polars of all members of a given equisingularity class. 

\begin{ejm} \rm Let $\Gamma=\langle 5,12\rangle$. The Normal Forms Theorem, together with the algorithm to compute normal forms in \cite{HH4}, gives us the complete classification of the curves in the equisingularity class determined by $\Gamma$ that we summarize in the following table:

{\footnotesize \begin{center}
\begin{tabular}{|l|l|}
\hline Normal Form &  $\Lambda_\ell\setminus \Gamma$ \\
\hline 
\hline 

1. \ $(t^5,t^{12})$ & $\emptyset$ \\ 
\hline 

2. \ $(t^5,t^{12}+t^{38})$ & $\{43\}$ \\ 
\hline 

3. \ $(t^5,t^{12}+t^{33})$& $\{38,43\}$\\ 
\hline 

4. \ $(t^5,t^{12}+t^{28})$& $\{33,38,43\}$\\ 
\hline 

5.  $(t^5,t^{12}+t^{26}+ct^{28}),\,\,\, c \neq 0 $&$\{31,38,43\}$\\ 
\hline 

6. \ $(t^5,t^{12}+t^{26}+ct^{33})$&$\{31,43\}$\\ 
\hline 

7. \ $(t^5,t^{12}+t^{23}+ct^{26})$& $\{28,33,38,43\}$\\ 
\hline 

8. \ $(t^5,t^{12}+t^{21}+ct^{23}+dt^{28})$&$\{26,31,38,43\}$\\ 
\hline 

9. \ $(t^5,t^{12}+t^{18}+ct^{21}+dt^{26})$&$\{23,28,33,38,43\}$\\ 
\hline 

10. \ $(t^5,t^{12}+t^{16}+ct^{18}+dt^{23})$&$\{21,26,31,33,38,43\}$\\ 
\hline 

11. \ $ (t^5,t^{12}+t^{14}+ct^{16}+dt^{18}+et^{23}), \,\, c\neq \frac{13}{12},\,\,d\neq \frac{4c^2-1}{3}$&$\{19,26,31,33,38,43\} $\\ 
\hline 

12. \ $(t^5,t^{12}+t^{14}+ct^{16}+(\frac{4c^2-1}{3})t^{18}+dt^{23}+et^{28}),\,\,c\neq \frac{13}{12}$&$\{19,26,31,38,43\}$\\ 
\hline 

13. \ $(t^5,t^{12}+t^{14}+\frac{13}{12}t^{16}+ct^{18}+dt^{21}),\,\, c\neq \frac{133}{108}$  &$\{19,28,31,33,38,43\}$\\ 
\hline 

14. \ $(t^5,t^{12}+t^{14}+\frac{13}{12}t^{16}+\frac{133}{108}t^{18}+ct^{21}+dt^{23}),\,\,d\neq \frac{34c}{11}$&$\{19,31,33,38,43\}$\\ 
\hline 

15. \ $(t^5,t^{12}+t^{14}+\frac{13}{12}t^{16}+\frac{133}{108}t^{18}+ct^{21}+\frac{34}{11}ct^{23}+dt^{28}),
$&$\{19,31,38,43\}$\\ 
\qquad $ d\neq\frac{81c^2}{32}+\frac{5225}{559872}$ &  \\ 
\hline 

16. \ $(t^{5},t^{12}+t^{14}+\frac{13}{12}t^{16}+\frac{133}{108}t^{18}+ct^{21}+\frac{34}{11}ct^{23}+
$&$\{19,31,43\}$\\   
\qquad  $(\frac{81c^2}{32}+\frac{5225}{559872})t^{28}+dt^{33})$& \\ 
\hline 

17. \,\, $(t^5, t^{12}+t^{13}-\frac 12 t^{14}+ct^{16}+dt^{21}+et^{26})$  & $\{18,23,28,33,38,43\} $ \\
\hline 

18. \ $(t^5,t^{12}+t^{13}+ct^{14}+dt^{16}+et^{21}),\,\,  c\neq -\frac 12$ & $ \{18,23,28,,31,33,38,43\} $ \\
\hline 

\end{tabular}\medskip

\footnotesize{Table 2.1: The normal forms of the equisingularity class of  $\langle 5, 12 \rangle$}
\end{center} }

Now, with the help of the Maple software we obtain the implicit equations of the curves given by the parametrization in each row of the table, then we exhibit their polars and analyze the equisingularity classes of these polars. In what follows, the symbols $u_1$, $u_2$, $u_3$ and $u_4$ represent units in $\C\{x,y\}$, with $u_i(0,0)=1$, for $1\leq i \leq 4$, not necessarily the same in all cases.\medskip 

\noindent {\bf 1.}  
$af_x+bf_y=5by^4-12ax^{11}.$ 

\noindent {\bf 2.} $af_x+bf_y=5by^4-50ax^9y^3-15bx^{10}y^2+100ax^{19}y-12ax^{11}+5bx^{20}-38ax^{37}.$ 

\noindent {\bf 3.} $af_x+bf_y=5by^4-45ax^8y^3-15bx^9y^2+ 90ax^{17}y-12ax^{11}+5bx^{18}-33ax^{32}.$ 

\noindent {\bf 4.} $af_x+bf_y=5by^4-40ax^7y^3-15bx^8y^2+80ax^{15}y-12ax^{11}+5bx^{16}-28ax^{27}.$ 

\noindent {\bf 5.} $af_x+bf_y=5by^4-40acx^7y^3-15bcu_1x^8y^2-10bu_2x^{10}y-12au_3x^{11}.$ 

\noindent {\bf 6.} $af_x+bf_y=5by^4-45acx^8y^3-(15bc+50a)u_1x^9y^2-10bu_2x^{10}y-12au_3x^{11}.$ 

\noindent {\bf 7.} $af_x+bf_y=5by^4-35ax^6y^3-15bu_1x^7y^2-10bcu_2x^{10}y-12au_3x^{11}.$

\noindent {\bf 8.} $af_x+bf_y=5by^4-35acu_1x^6y^3-15bcu_2x^7y^2-10bu_3x^9y-12au_4x^{11}.$ 

\noindent {\bf 9.} $af_x+bf_y=5by^4-30ax^5y^3-15bu_1x^6y^2-10bcu_2x^9y-12au_3x^{11}.$ \smallskip

The polar of any one of the curves in the families $(1)-(9)$ has Newton polygon with only one side $L=[(0,4);(11,0)]$ that supports only its extremal points associated to the monomials  $5by^4$ and $-12ax^{11}$. This implies that all (general) polars are Newton non-degenerate, so their Newton polygons determine their equisisingularity classes, which in this case is given by only one branch with semigroup $\langle 4,11\rangle$. \medskip

\noindent {\bf 10.} $af_x+bf_y=5by^4-30acu_1x^5y^3-15bcu_2x^6y^2-10bu_3x^8y-12au_4x^{11}.$\smallskip

In this case, the Newton polygon of the polar has two sides:\smallskip

\noindent $L_1=[(0,4);(8,1)]$, that supports only its extremal points associated to the monomials  $5by^4$ and $-10bx^8y$. \smallskip

\noindent $L_2=[(8,1);(11,0)]$, that supports only its extremal points associated to the monomials $-10bx^8y$ and $-12ax^{11}$. \smallskip

Again, the polar of any curve belonging to this family is Newton non-degenerate with two branches: $g_1$ with semigroup $\langle 3,8 \rangle$ and $g_2$ smooth such that ${\rm I}(g_1,g_2)=8$.\medskip

\noindent {\bf 11.} 

$\begin{array}{rl} af_x+bf_y&=
5by^4-30a(c+d)u_1x^5y^3-15b(c+d)u_2x^6y^2 \\
& \ -10b(1+c)u_3x^8y-5bu_4x^{10}.\end{array}$ \smallskip

\noindent {\bf 12.} 

$\begin{array}{rl}  af_x+bf_y&=5by^4+10a(1-3c-4c^2)u_1x^5y^3+5b(1-3c-4c^2)u_2x^6y^2\\ &  \ -10b(1+c)u_3x^8y-5bu_4x^{10}.
\end{array}$\smallskip

\noindent {\bf 13.} 

$\begin{array}{rl} af_x+bf_y&=5by^4-\frac 52a(13+12c) u_1x^5y^3-\frac 54 b(13+12c)u_2x^6y^2 \\ & \ -\frac{125}{6}bu_3x^8y -5bu_4x^{10}.
\end{array}$ \smallskip

\noindent {\bf 14.}  $af_x+bf_y=5by^4-\frac{625}{9}au_1x^5y^3-\frac{625}{18}bu_2x^6y^2-\frac{125}{6}bu_3x^8y -5bu_4x^{10}.$ \smallskip

\noindent {\bf 15.}  $af_x+bf_y=5by^4-\frac{625}{9}au_1x^5y^3-\frac{625}{18}bu_2x^6y^2-\frac{125}{6}bu_3x^8y-5bu_4x^{10}.$ \smallskip

\noindent {\bf 16.}  $af_x+bf_y=5by^4-\frac{625}{9}au_1x^5y^3-\frac{625}{18}bu_2x^6y^2-\frac{125}{6}bu_3x^8y-5bu_4x^{10}.$ \smallskip

\noindent {\bf 17.}
$af_x+bf_y=5by^4-25au_1x^4y^3-15bu_2x^5y^2+5b\big(\frac{1-4c}{2}\big)u_3x^8y+\frac{15}{2}bu_4x^{10}.$\medskip

The Newton polygon of the polars of any member of the families $(11)-(17)$ has only one side  $L=[(0,4);(10,0)]$, that supports just its extremal points associated to the monomials $y^4$ and $x^{10}$ with some non-zero coefficients that do not depend upon the parameters $c,d$ and $e$. Therefore, these polars are Newton non-degenerate, so they have two branches with semigroup $\langle 2,5\rangle$ that intersect with multiplicity $10$. \medskip

\noindent {\bf 18.}

$\begin{array}{ll}
af_x+bf_y =& 5by^4-25au_1x^4y^3 - 15bu_2x^5y^2\\
           &-10[b(c^2+c+d)x^8+(b(d^2+e)+5a(c-1))u_3x^9]y \\
           & +5b(1-c)x^{10}+(-5b(c^3+dc+d)-12au_4x^{11}.
\end{array}$
\smallskip

This is the only stratum in which the equisingularity class of the polars will depend upon the parameters in ${\mathcal E}_\ell$ to which it belongs. \medskip

\noindent (i) If $c\neq 1$, then the Newton polygon of the polar has the only side $L=[(0,4);(10,0)]$, that supports the points associated to the monomials $5by^4$, $-15bx^5y^2$ and $-5b(c-1)x^{10}$. In this case, $p_L(z)=z^4-3z^2-(c-1)$ whose discriminant is $-16(c-1)(5+4c)^2$. So, the polar is Newton non-degenerate if and only if $c\neq -\frac{5}{4}$. In this case, the polar has two branches with semigroup $\langle 2,5\rangle$ that intersect with multiplicity $10$.\medskip

If $c=-\frac{5}{4}$, the parametrization of the polar is given by 
$$x=\frac{3^{20}}{2^{21}}(16d+5)^2t^4,\,\,\, y=\frac{3^{38}}{2^{53}}(16d+5)^5t^{10}+ \frac{3^{30}}{2^{60}}(16d+5)^6t^{11}+\cdots.$$

Therefore, when $d \neq -\frac{5}{16}$, the members of the family ${\mathcal E}_{18}$ for which $c=-\frac{5}{4}$ have irreducible polars of genus $2$ with semigroup $\langle 4,10,21\rangle$,  and when $d=-\frac{5}{16}$, they have  
polars with two branches with parametrizations
\[
x_i=\frac 32t^2, \ \ y_i=\frac{27}{8}t^5+(-1)^i\frac{27}{640b}(256ab-125b)^{\frac 12}t^6+\cdots, \ \ i=1,2,
\]
that is, branches with semigroup  $\langle 2,5\rangle$ and with intersection multiplicity $11$.\medskip

\noindent (ii) If $c=1$, then the polar is given by
$$\begin{array}{l} 
af_x+bf_y=\\ \\
5by^4-25ax^4u_1y^3-15bx^5u_2y^2-10b(d+2)x^8u_3y-(12a+5b+10bd)x^{11}u_4.
\end{array}$$

\noindent An easy computation shows that its Newton polygon has two sides:\smallskip

\noindent $L_1=[(0,4);(5,2)]$, that supports only its extremal points associated to the monomials $5by^4$ and $-15bx^5y^2$; and\smallskip

\noindent $L_2=[(5,2);(11,0)]$, that supports only its extremal points associated to the monomials $-15bx^5y^2$ and $-(12a+5b+10bd)x^{11}$.\smallskip

Therefore, any curve in this family with $c=1$ has, for general values of $a$ and $b$, a Newton non-degenerate polar with a branch $p$ with semigroup $\langle 2,5\rangle$ and two non-singular branches $g_1$ and $g_2$ such that ${\rm I}(p,g_i)=5$ and ${\rm I}(g_1,g_2)=3$.\smallskip
\end{ejm}

\noindent {\bf Remark} \ The stratum ${\mathcal E}_{18}$ gives us an example in which the equisingularity class of the general polar of its members is not constant. It also gives us a somewhat unexpected example of a family of curves of genus $1$ such that its general member has a general polar of genus $2$.\medskip

What is remarkable is that the analytic classification of the branches in this equisingularity class allowed us to describe the equisingularity classes of all general polars of its members.\medskip

Although, as we saw in the above example, the topological type of the polar may be not constant in a given stratum ${\mathcal E}_\ell$, it is constant in an open dense set of each irreducible component of the stratum, as we will show in general in the sequel.

In fact, for the stratum associated to $\Lambda=\Gamma\setminus\{0\}$, the result follows easily. Let us consider a normal form in an equisingularity class para\-me\-tri\-zed by ${\mathcal E}_{\ell}$ associated to a set of values of differentials $\Lambda_\ell \neq \Gamma\setminus\{0\}$. Putting $v_0=n$ and $v_1=m$, from the Normal Forms Theorem, we have
$$x=t^n,  \ \ y = t^m+ t^\lambda+\sum_{{\stackrel{i> \lambda}{i \notin \Lambda_\ell -n}}} c_i t^i.$$

The implicit equations of these curves are given by Weierstrass polynomials
$$f=y^{n}+a_2(x)y^{n-2}+a_3(x)y^{n-3}+\cdots +a_{n-1}(x)y+a_n(x),$$
where the coefficients $a_j(x)$ are polynomials in the variables $c_i$ and such that 
${\rm ord}_x(a_j(x)) > j$ and ${\rm ord}_xa_{n}(x)=m$.\medskip

Therefore, the polars of the curves in ${\mathcal E}_\ell$ are given by the family 
$$
\begin{array}{rcl}
P(f)&=&af_x+bf_y \\
    &=& bny^{n-1}+aa_2'(x)y^{n-2}+\big(ba_2(x)(n-2)+aa_3'(x) \big)y^{n-3} +\cdots \\
    &&  +\big (b(n-2)a_{n-2}(x)+aa_{n-1}'(x)\big)y+ba_{n-1}(x)+aa_{n}'{x}.
\end{array}
$$

We will now show that in a dense open Zariski set in any irreducible component of ${\mathcal E}_\ell$ the value of the Milnor number of $P(f)$ is constant. \medskip

From the equation of $P(f)$ we have that
$$ 
\begin{array}{rcl}
P(f)_x&=&aa_2''(x)y^{n-2}+\big( aa_3''(x)+b(n-2)a_2'(x) \big )y^{n-3}+ \cdots +\\
& & \big(aa_{n-1}''(x)+2ba_{n-2}'(x)\big )y+aa_n''(x)+ba_{n-1}'(x), \quad \text{and} \\ \\
P(f)_y&=&bn(n-1)y^{n-2}+a(n-2)a_2'(x)y^{n-3}+\cdots + \\
& & \big( aa_{n-1}'(x)+b(n-2)a_{n-2}(x) \big).
\end{array}$$

\noindent Therefore, one has that $P(f)_y$ is a constant times a Weierstrass polynomial in $y$ and $P(f)_x \in \C\{x\}[y]$, hence their intersection multiplicity, which is the Milnor number of $P(f)$, is the order in $x$ of their resultant $R_y$ in $y$. Because $R_y\not \equiv 0$ since the generic polar of $f$ is reduced, we have for every irreducible component ${\mathcal E}_{\ell,j}$ of ${\mathcal E}_\ell$ that
$$R_y \big(P(f)_x,P(f)_y \big)=A_jx^{\nu_j} +{higher \ order \ terms},$$
where $A_j$ is a non-zero polynomial in $a$, $b$ and the $c_i$ (the coefficients in the normal forms)  and homogeneous in $a$ and $b$. So, there exists a Zariski open set in ${\mathcal E}_{\ell,j}$, where at each point this polynomial $A_j$ in $a$ and $b$ is not identically zero, hence the Milnor number of the general polar of the corresponding curve is constant $(=\nu_j)$ in this open set. From the  L\^e-Ramanujan Theorem \cite{LR}, we obtain the following result:
\begin{teor}
The equisingularity class of the polar of curves in ${\mathcal E}_\ell$ is constant in an open dense Zariski subset of any of its irreducible components.
\end{teor}

\section{Polars of branches up to multiplicity four}

We will now give a detailed description of the equisingularity classes of the polars of branches of multiplicity less or equal than four. This will be carried out by using the classification done by the first two authors in \cite{HH3}. Observe that the polar of a branch of multiplicity $2$ is a smooth branch, so we have only to treat the cases of multiplicities three and four.

\subsection{Multiplicity three} 
For multiplicity three curves, there is only one analytic representative in each stratum which is determined by Zariski's $\lambda$ invariant, as shown in the table below:
\begin{center}
\begin{tabular}{|l|}
\hline $\,\,\,\,\Gamma= \langle 3,\beta\rangle;   \ \ \,\,\,\,\,\,\, \beta=3q +\varepsilon$, 
\ \ $\varepsilon=1,\; 2$ \\ \hline
\hline $x=t^3,\,\,\,y=t^\beta$ \\
\hline $x=t^3,\,\,\,y=t^\beta+t^{\beta+\varepsilon +3k},\,\,0 \leq k \leq q-2$\\
\hline
\end{tabular}
\end{center}
For the case of the monomial curve  $x=t^3,\,\,y=t^\beta$, 
we have that the polar curve has $d=\tx{gcd}(2,\beta-1)$ branches.
When $d=1$, the branch has semigroup $\langle 2,\beta-1\rangle$ and when $d=2$, the two branches are smooth and their intersection multiplicity is $\frac{\beta-1}{2}$.\medskip

In the case of the second row of the above table, the implicit equation of the curve is
$f=y^3-3x^{2q+k+\epsilon}y-x^{\beta}-x^{\beta+\epsilon+3k}$ and the generic polar curve is
$$ af_x+bf_y=2by^2-3a(2q+k+\epsilon)x^{2q+k+\epsilon}y-3bux^{2q+k+\epsilon},
$$
where $u$ is a unit. After a direct computation, we see that the equisingularity class of the polar may be described by the following table:
{\footnotesize \begin{center}
\begin{tabular}{|l|l|}
\hline $2q+k+\epsilon=2I+1$ & One branch with semigroup $\langle 2,2q+k+\epsilon\rangle$.\\ 
\hline $2q+k+\epsilon=2I$ & Two smooth branches with intersection multiplicity $I$.\\  \hline
\end{tabular}
\end{center}}

\subsection{Multiplicity four, genus one}
A curve of multiplicity $4$ may have genus one or two. For the genus one case, we have the following normal forms:
{\footnotesize{\begin{center}
\begin{tabular}{|l|l|}
\hline Normal form  &  $\Lambda\setminus \langle 4,m \rangle$ \\
\hline 1. \ $y(t)=t^{m}$ &  $\emptyset$\\
\hline
2. \ $y(t)=t^m+t^{3m-4j}+a_1t^{2m-4(j-[\frac{m}{4}]-1)}+ \cdots $ & $ \{3m-4s; 1\leq s \leq j-1\}$ \\
\qquad $+ a_{j -[\frac{m}{4}]-2}t^{2m-8}$, \,\,\,\,\,\,\,\,\, $2 \leq j \leq [\frac{m}{2}]$& \\
\hline
3. \ $y(t)=t^{m}+t^{2m-4j}+a_kt^{3m-(4[\frac{m}{4}]+j+1-k)}+ \cdots$ &  $\{ 2m-4s; 1 \leq s \leq j-1\} \ \cup $ \\
\qquad  $+ a_{j-[\frac{m}{4}]-2}t^{3m-4([\frac{m}{4}]+3-k)}$ & $\{ 3m-4s;1\leq s\leq[\frac{m}{4}]+1-k \}$\\
\qquad $a_k \neq 0, \,\, \,\,\,\, \,\,\,  2\leq j \leq [\frac{m}{4}], \,\,\,\,   1\leq k \leq [\frac{m}{4}]-j$ &\\
\hline
4. \ $y(t)=t^{m}+t^{2m-4j}+a_{[\frac{m}{4}]-j+1}t^{3m-8j}$& $\{2m-4s; 1 \leq s \leq j-1\} \ \cup$\\
\qquad $ +a_{[\frac{m}{4}]-j+2}t^{3m-4(2j-1)}+ \cdots + a_{[\frac{m}{4}]-1}t^{3m-4(j+2)}$,&   $\{3m-4s;1 \leq s \leq j \}$\\
$\,\,\,\,\,\,\,\,\,a_{[\frac{m}{4}]-j+1} \neq \frac{3m-4j}{2m}, \,\,\,\,\,\,\,\,\,  2 \leq j \leq [\frac{m}{4}]$& \\
\hline
5. \ $y(t)=t^{m}+t^{2m-4j}+\frac{3m-4j}{2m}t^{3m-8j}+$&
$\{2m-4s;1\leq s\leq j-1\} \ \cup$\\
\qquad $a_{[\frac{m}{4}]-j+2} t^{3m-4(2j-1)}+\cdots+ a_{[\frac{m}{4}]}t^{3m-(j+1)}$, & $\{ 3m-4s;1\leq s\leq j-1\}$ \\
\qquad $2\leq j \leq [\frac{m}{4}]$ &   
\\ \hline
\end{tabular}
\end{center} 
\vspace{.3cm}
\centerline{Table 3.1: Normal forms for multiplicity four and genus one}} }\medskip

\noindent {\bf First Normal Form} (monomial curves)\medskip

In this case, the equation of the curve is $y^4-x^m=0$, so its polar is $4by^3-amx^{m-1}$, that has $d={\rm gcd}(3,m-1)$ branches. If $d=1$, the branch has semigroup $\langle 3,m-1\rangle$ and when $d=3$, the three branches are smooth with mutual intersection multiplicity equal to $\frac{(m-1)}{3}$.
\medskip

\noindent {\bf Second Normal Form}\medskip

This is the more complicated case. The implicit equation of the curve is
$$f=y^{4}- S_1(x)y^{3}+S_2(x)y^{2}-S_3(x)y + S_4(x)=0,$$
where $S_r(x)$ is the $r$-th symmetric polynomial computed in $y(\varepsilon^lt)$, $l=0,1,2,3$, with $\varepsilon$ a primitive fourth root of $1$ and where we have replaced $t^4$ by $x$.

From the definition of $y(t)$, it is clear that $S_1=0$. To determine the Newton polygon of the polar, it is sufficient to consider in the polynomial $S_r(x)$, $2\leq r\leq 3$, the monomial which determines its multiplicity.\\

\noindent {\bf (I)} \ We first consider the case $a_1=a_2= \ldots =a_{j-[\frac{m}{4}]-2}=0$.
\medskip

For each fixed $j$, we have $f=y^4-4x^{m-j}y^2-x^m+2x^{2m-2j}-x^{3m-4j}$. 
Therefore,
$$af_x+bf_y=4by^3-4a(m-j)x^{m-j-1}y^2-8bx^{m-j}y-amx^{m-1}u,$$
where $u\in\mathbb{C}\{x\}$ with $u(0)=1$.\medskip

We have the following cases:\medskip

\noindent {\bf i.} \ Case $\frac{2}{m-j} < \frac{1}{j-1}$.\medskip

In this case, the Newton polygon of the polar has only one side $L$ containing only its end points $(0,3)$ and $(m-1,0)$, associated to monomials of the polar. 
The polynomial associated to the Newton polygon is $p_L(z)=4bz^3-am$.
Then for $a$ and $b$ generic, $p_L(z)$ has three distinct roots 
 $\{z_1,z_2,z_3\}$. Therefore the polar has:\medskip

\noindent a) One branch with semigroup $\langle 3,m-1\rangle$, if $\tx{gcd}(3,m-1)=1$.\medskip

\noindent b) Three smooth branches with parametrizations: $(t, z_it^{\frac{m-1}{3}}+\cdots )$,
and mutual intersection numbers $\frac{m-1}{3}$, if
$\tx{gcd}(3,m-1)=3$.
\medskip

\noindent {\bf ii.} \ Case $\frac{2}{m-j} > \frac{1}{j-1}$.\medskip

In this case, the Newton polygon of the polar has two sides $L_1$ and $L_2$,
each one with only its end points associated to monomials of the polar.
The associated polynomials are $p_{L_1}(z)=4bz^2-8b$, and $p_{L_2}(z)=-8bz-am$. Then, we have that \medskip

\noindent a) Associated to $L_1$ there is one branch $p_1$ with semigroup $\langle 2,m-j\rangle$ and parametrization
$x=t^2\,\,\,y=\sqrt{2}t^{m-j}+\cdots$, if $\tx{gcd}(2,m-j)=1$; or two smooths  branches  $g_1,g_2$ with parametrizations
$x_1=t,\,\,\,y_1=\sqrt{2}t^{\frac{m-j}{2}}+\cdots $ and
$x_2=t,\,\,\,y_2=-\sqrt{2}t^{\frac{m-j}{2}}+\cdots $, if $\tx{gcd}(2,m-j)=2$.\medskip

\noindent b) Associated to $L_2$, there is one smooth branch $p_2$ with parametrization
 $x=t,\,\,y=-\frac{am}{8b}t^{j-1}+\cdots$. \smallskip

Finally, we have that ${\rm I}(p_1,p_2)=m-j$ and 
${\rm I}(g_i,p_2)={\rm I}(g_1,g_2)=\frac{m-j}{2}$.\medskip

\noindent {\bf iii.} \ Case $\frac{2}{m-j} = \frac{1}{j-1}$.\medskip

Since $j>2$, because otherwise $m=4$, which is not allowed, the Newton polygon of the polar has only one side $L$ with tree points and the polynomial associated to $L$ is $p_L(z)=4bz^3-8bz-am$. Therefore, for $a$ and $b$ generic, the polynomial 
$p_L(z)$ has three distinct roots $\{z_1,z_2,z_3\}$ and as, in this case,  $\tx{gcd}(3,m-1)=3$, then associated to $L$ we have three smooth branches with parametrizations 
$(t, z_it^{j-1}+\cdots)$, $i=1,2,3$, and mutual intersection numbers $j-1$.\medskip

\noindent {\bf (II)} \ Now we consider the case where some of the $a_i$'s is non-zero. Set 
$k= \tx{min}\{i ; \ a_i \neq 0\}$. \medskip

\noindent After a computation we get 
$$f=y^4-(4x^{m-j}+2a_k^2x^{m-2(j-\left [\frac{m}{4} \right ]-k)}u_1)y^2-4a_kx^{m-(j-\left [\frac{m}{4} \right ]-k)}u_2y-x^mu_3,$$ where $u_i\in\mathbb{C}\{x\}$ with $u_i(0)=1$ for $i=1,2,3$. Hence, to determine the Newton polygon of the polar 
$af_x+bf_y$, it is sufficient to consider the polynomial
$$4by^{3}-4a(m-j)x^{m-j-1}y^{2}-8bx^{m-j}y-4ba_kx^{m-j+[\frac{m}{4}]+k}.$$

We now split the analysis of this case into several sub-cases.\medskip

\noindent {\bf i.} \  Case $\frac{2}{m-j} < \frac{1}{[\frac{m}{4}]+k}$.

The Newton polygon of the polar has just one side $L$, containing only the points $(0,3)$ 
and $(m-j+[\frac{m}{4}]+k,0)$.

Since the polynomial $p_L(z)=4bz^3-4ba_k$ has three distinct roots $\{z_1,z_2,z_3\}$, it follows that the polar has:\smallskip

\noindent a) Only one branch, if $\tx{gcd}(3,m-j+[\frac{m}{4}]+k)=1$, with semigroup $\langle 3,m-j+[\frac{m}{4}]+k\rangle$.
\medskip

\noindent b) Three  smooth branches, if $\tx{gcd}(3,m-j+[\frac{m}{4}]+k)=3$, with parameterizations
 $$x_i=t,\,\,\;y_i=z_it^{\frac{m-j+[\frac{m}{4}]+k}{3}}+\cdots, \quad i\in \{1,2,3\},
$$
and mutual intersection numbers $\frac{m-j+[\frac{m}{4}]+k}{3}$.\medskip

\noindent {\bf ii.} \ Case $\frac{2}{m-j} > \frac{1}{[\frac{m}{4}]+k}$.

In this case, the Newton polygon of the polar has two sides $L_1=[(0,3);(m-j,1)]$ and
 $L_2=[(m-j,1);(m-j+[\frac{m}{4}]+k,0)]$, with on each side only the extreme points correspond to monomials of the polar.

Considering the polynomials associated to these sides, 
$p_{L_1}(z)=4bz^2-8b$ and $p_{L_2}(z)=-8bz-4ba_k$; and defining $d=\tx{gcd}(2,m-j)$, we have that\smallskip

\noindent a)  Associated to the side $L_1$, we have a branch $p_1$ with semigroup $\langle 2,m-j\rangle$ and parametrization $x=t^2,\,\,\,\,y=\sqrt{2}t^{m-j}+\cdots$, if $d=1$; and two  smooth branches $g_1,g_2$, with parametrizations
$x_i=t,\,\,\, y_i= (-1)^{i-1} \sqrt{2}t^{\frac{m-j}{2}}+\cdots $, $i=1,2$, if $d=2$.\smallskip
 
\noindent b) Associated to the side $L_2$, we have a smooth branch $p_2$, with parametrization
 $x=t,\,\,\, y=-\frac{a_k}{2}t^{[\frac{m}{4}]+k}+\cdots$.\smallskip
 
Finally, one has ${\rm I}(p_1,p_2)=m-j$ \,\, and \,\,
${\rm I}(g_i,p_2)={\rm I}(g_1,g_2)=\frac{m-j}{2}$.\medskip

\noindent {\bf iii.} \ Case $\frac{2}{m-j}=\frac{1}{[\frac{m}{4}]+k}$. \medskip

In this case, the Newton polygon of the polar has a unique side $L$ containing 
the three points $(0,3),(m-j,1)$ and $(m-j+[\frac{m}{4}]+k,0)$.
whose associated polynomial is $$p_L(z)=4bz^3-8bz-4a_k.$$

When $a_k\neq \frac{4\sqrt{6}}{9}(-1)^{\alpha};\ \alpha=0,1$, because of the condition $\frac{2}{m-j}=\frac{1}{[\frac{m}{4}]+k}$, it is easy to verify that the polar is Newton non-degenerate. In this case, the polynomial $p_L(z)$ has three distinct roots $\{z_1,z_2,z_3\}$, then the polar has three smooth branches with parametrizations $x=t$ and $y_i= z_it^{[\frac{m}{4}]+k}+\cdots$, \ $i=1,2,3$, with mutual intersection numbers equal to $[\frac{m}{4}]+k=\frac{m-j}{2}$.\medskip

Now we suppose that $a_k= \frac{4\sqrt{6}}{9}(-1)^{\alpha};\ \alpha=0,1$.\medskip

In this case, the roots of $p_L(z)$ are $\frac{\sqrt{6}}{3}(-1)^{\alpha+1},\frac{\sqrt{6}}{3}(-1)^{\alpha+1}$ and $2\frac{\sqrt{6}}{3}(-1)^{\alpha}$. The polar will have a smooth branch $f_1$ corresponding to the simple root of $p_L(z)$ and branches $g_i$ corresponding to the double root.

We may suppose that the roots of $p_L(z)$ are $\frac{\sqrt 6}{3}, \frac{\sqrt 6}{3}$ and $-2\frac{\sqrt 6}{3}$, since the other case is analogous. \medskip

\noindent a) If for all $l>0$ one has $a_{k+l}=0$, then a simple analysis shows that the polar has a smooth branch and a branch with semigroup $\langle 2, 2m-3j \rangle$ with intersection number $m-j$.\medskip

\noindent b) Suppose that there exists $l>0$ such that $a_{k+l}\neq 0$. We denote the least such $l$ by $s$. In this case, we will need in our analysis to consider more terms of $f$, which now reads  as
\[
\begin{array}{rcl}
f &=& y^4+(-4x^{m-j}-2a_k^2 x^{m-2(j-[ \frac m4 ]-k)}-4a_k a_{k+s}x^{m-2(j-[ \frac m4 ]-k)+s}\\
  & & -2a_{k+s}^2x^{m-2(j-[ \frac m4 ]-k-s)}+\cdots )y^2-(4a_kx^{m-j+[ \frac m4 ]+k} \\
&& +4a_{k+s}x^{m-j+[ \frac m4 ]+k+s} +4a_kx^{2m-2j-(j- [ \frac m4 ]-k)}\\
 && +4a_{k+s}x^{2m-2j-(j-[ \frac m4 ]-k-s)}+\cdots )y-ux^m,
\end{array}
\]
where $u\in \C\{x\}$ with $u(0)=1$.

Now, in order to apply the Newton-Puiseux algorithm to the general polar of $f$ at the double root of $p_L(z)$,  we have to split our analysis in several subcases.\medskip

\noindent  b.1) \ $m-2j>s$.\medskip

\noindent b.1.1) \ $s$ odd. Associated to the double root there is a branch $g_1$ given by
\[
x=t^2,\ y=-\frac{\sqrt{6}}{3}t^{m-j}+\frac{\sqrt{a_{k+s}}}{\sqrt[4]{6}}t^{m-j+s}+\ldots.
\]
In this case, the polar has a smooth branch $f_1$ and a branch $g_1$ with semigroup $\langle 2, m-j+s \rangle$ such that ${\rm I}(f_1,g_1)=m-j$. \medskip

\noindent b.1.2) \ $s$ even. The polar splits into three smooth factors $f_1, g_1$ and $g_2$, such that ${\rm I}(f_1,g_i)=[\frac m4]+k$ and ${\rm I}(g_1,g_2)=[\frac m4]+k+\frac s2$.\medskip

\noindent  b.2) \ $m-2j<s$. In this case, the polar has the smooth branch $f_1$ and a branch  branch $g_1$ associated to the double root with semigroup $\langle 2, 2m-3j\rangle$, such that  ${\rm I}(f_1,g_1)=m-j$.\medskip

\noindent  b.3) \ $m-2j=s$.\medskip

\noindent  b.3.1) \ If $a_{k+s}\neq \frac{4\sqrt{6}}{81}(-1)^{\alpha+1}$, we have, associated to the double root, a branch $g_1$ with semigroup $\langle 2, m-j+s \rangle$. So, the polar has the smooth branch $f_1$ and the above branch $g_1$ such that ${\rm I}(f_1,g_1)=m-j$.\medskip

\noindent  b.3.2) \ If $a_{k+s}= \frac{4\sqrt{6}}{81}(-1)^{\alpha+1}$, we have, associated to the double root, two smooth branches $g_1$ and $g_2$ such that ${\rm I}(f_1,g_i)=\frac{m-j}{2}$, $i=1,2$, and ${\rm I}(g_1,g_2)=\frac{m-j}{2}+s$.\medskip

The following table summarizes the above analysis for the second normal form.

{\footnotesize \begin{center}
\begin{tabular}{|l|}
\hline  $y=t^m+t^{3m-4j}+a_1t^{2m-4(j-[\frac{m}{4}]-1)}+ \cdots + a_{j - [\frac{m}{4}]-2}t^{2m-8};\,\, 2 \leq j \leq [\frac{m}{2}]$\\
\hline
\hline 
\centerline{  \ $a_1=a_2=\cdots= a_{j-[\frac{m}{4}]-2}=0$} \\
\hline
\begin{tabular}{l|l}
              & The polar has one branch with semigroup  $\langle 3,m-1\rangle$, \\
$\frac{2}{m-j}<\frac{1}{j-1}$ &if $\tx{gcd}(3,m-1)=1$; \ otherwise it has three smooth branches \\
& $p_1,p_2,p_3$, with ${\rm I}(p_i,p_r)=\frac{m-1}{3}$.\\ \hline
& The polar has one branch $p_1$ with semigroup  $\langle 2,m-j\rangle$ and \\
$\frac{2}{m-j}>\frac{1}{j-1}$& one smooth branch $p_2$, with ${\rm I}(p_1,p_2)=m-j$, if $\tx{gcd}(2,m-j)=1$;\\
& otherwise it has three smooth branches $p_1,p_2,p_3$, with ${\rm I}(p_i,p_r)=\frac{m-j}{2}$. \\ 
\hline 
$\frac{2}{m-j}=\frac{1}{j-1}$ & The polar has three smooth branches   $p_1,p_2,p_3$, with ${\rm I}(p_i,p_r)=j-1$. \\ 
\hline 
\end{tabular}\\
\centerline{ $\exists\; i; \ a_i\neq 0, 
\,\,\,\,\,\,\,\,\,\,\,\,\,\ \ \ k=\tx{min}\{i; \ a_i\neq0 \}$ } \\
\hline
\begin{tabular}{l|l}
& The polar has one branch with semigroup  $\langle 3,m-j+[\frac{m}{4}]+k\rangle$,\\
 				$\frac{2}{m-j}< \frac{1}{[\frac{m}{4}]+k}$  &if ${\rm gcd}(3,m-j+[\frac{m}{4}]+k)=1$; otherwise it has three smooth branches\\
 				&  $p_1,p_2,p_3$ with ${\rm I}(p_i,p_r)=\frac{m-j+[\frac{m}{4}]+k}{3}$.\\ \hline
 				& The polar has a branch $p_1$, with semigroup $\langle 2,m-j\rangle$\\
 				$\frac{2}{m-j}> \frac{1}{[\frac{m}{4}]+k}$ & and a smooth branch $p_2$, with ${\rm I}(p_1,p_2)=m-j$,\,\, if ${\rm gcd}(2,m-j)=1$;\\
 				& otherwise it has three smooth branches $p_1,p_2,p_3$, with ${\rm I}(p_i,p_r)=\frac{m-j}{2}$. \\ \hline

& For $a_k\neq \frac{4\sqrt{6}}{9}(-1)^{\alpha};\ \alpha=0,1$, the polar has three smooth branches  \\ 
 				& \quad $p_1,p_2,p_3$,  with  ${\rm I}(p_i,p_r)=\frac{m-j}{2}$.\\ 
 				&For $a_k= \frac{4\sqrt{6}}{9}(-1)^{\alpha};\ \alpha=0,1$:\\
 				&a) If $a_{k+l}=0$, $\forall l>0$, then the polar 
                 has a smooth branch $f_1$ \\
&and a branch $g_1$ with semigroup  $\langle 2,2m-3j\rangle$ with ${\rm I}(f_1,g_1)=m-j$.\\
 				&b) There exists $s>0$ such that $a_{k+s}\neq 0$ (let $s$ be minimum).\\	
 				& b.1) \ $m-2j>s$.\\ 
 				& b.1.1) $s$ odd. The polar has a smooth branch $f_1$                  and a branch $g_1$\\
 				&  with semigroup $\langle 2, m-j+s \rangle$ with ${\rm I}(f_1,g_1)=m-j$.\\
 				& b.1.2)  $s$ even. The polar has three  smooth  branches $f_1,g_1,g_2$ with \\
 				& ${\rm I}(f_1,g_i)=\frac{m-j}{2}$ and ${\rm I}(g_1,g_2)=\frac{m-j+s}{2}$.\\
 				$\frac{2}{m-j}=\frac{1}{[\frac{m}{4}]+k}$ & b.2) \ $m-2j<s$.  The polar has a smooth branch $f_1$ and a branch  $g_1$ \\ 
 				&  with semigroup  $\langle 2, 2m-3j\rangle$ with ${\rm I}(f_1,g_1)=m-j$.\\
 				& b.3) \ $m-2j=s$.\\ 
 				& b.3.1) $a_{k+s}\neq \frac{4\sqrt{6}}{81}(-1)^{\alpha+1}$. The polar has a smooth branch $f_1$ and a \\
 				& branch $g_1$ with semigroup $\langle 2 , 2m-3j\rangle$ with ${\rm I}(f_1,g_1)=m-j$.\\
 				& b.3.2)  $a_{k+s} = \frac{4\sqrt{6}}{81}(-1)^{\alpha+1}$. The polar has three smooth branches \\ 
 				 				& $f_1,g_1,g_2$  such that ${\rm I}(f_1,g_i)=\frac{m-j}{2}$ and ${\rm I}(g_1,g_2)=\frac{m-j}{2}+s$.\\
 				\hline
 \end{tabular}
\end{tabular}
\end{center} \medskip
\centerline{Table 3.2: The polars for curves in the second Normal Form}}\bigskip

\noindent {\bf Third to Fifth Normal Forms}\medskip

These are simple to analyze and all give the same result, summarized in the following table:

{\footnotesize\begin{center}
\begin{tabular}{|l|l|}
\hline
$\tx{gcd}(3,m-j)=1$ & One branch with semigroup $\langle 3,m-j\rangle$. \\ \hline
 $\tx{gcd}(3,m-j)=3$ & Three smooth branches with mutual intersection numbers $\frac{m-j}{3}$.\\ \hline
\end{tabular}
\end{center}
\centerline{Table 3.3: The polars for curves in the third to fifth Normal Forms}}\bigskip

\subsection{Multiplicity four and genus two}

The classification of multiplicity $4$ and genus $2$ branches is given in the table below, extracted from \cite{HH2}.

\begin{center}
\begin{tabular}{|l|l|}
\hline Normal form  & $ \Lambda\setminus \langle 4,v_1,v_2\rangle$\\
\hline 
$y(t)=t^{v_1}+t^{v_2 -v_1}+a_1t^{v_2-4[\frac{v_1}{4}]}$ & $v_2+v_1-4s;$ \\
\qquad $+a_2t^{v_2-4([\frac{v_1}{4}]-1)}+\cdots +a_{[\frac{v_1}{4}]-1}t^{v_2-8}$, & $1\leq s\leq \frac{v_1}{2}+1$ \\ 
\hline
\end{tabular}
\end{center}
\vspace{.3cm}
\centerline{\footnotesize{Table 3.4: Multiplicity four and genus two}}\medskip

Since ${\rm gcd}(4,v_1,v_2)=1$ and $v_2>2v_1$, we may write $v_1=2k_1$ and $v_2=4k_1+d$, where $k_1$ and $d$ are odd numbers. It is easy to verify that $2v_2-v_1$ is multiple of $4$, let us write $k_2=\frac{2v_2-v_1}{4}$. 

Writing an implicit equation for the curve, we have
\[
f=y^4+(-2x^{k_1}+\cdots)y^2+(-4x^{k_2}+\cdots)y+x^{v_1}u,
\]
where $u\in\mathbb{C}\{x\}$ and $u(0)=1$.

The Newton polygon of the polar $af_x+bf_y$ is determined by the polynomial
\[
4by^3-2k_1ax^{k_1-1}y^2-4bx^{k_1}y-4bx^{k_2}+v_1ax^{v_1-1}.
\]

This polygon has always two sides $L_1$ and $L_2$ whose positions depend upon the relationship among $k_2$ and $v_1-1$. In all cases, the polar will have a component $g_1$ associated to side $L_1$ given by the following parametrization
$$x_1=t^2,\,\,\, y_1=t^{k_1}+\cdots.$$

Associated to the side $L_2$, the polar has a branch $g_2$ parametrized by
$$
x_2=t, \ \ y_2=\left\{ \begin{array}{ll}
\frac{ak_1}{2b}t^{v_1-k_1-1}+\cdots, & \mbox{if}\ v_1-1<k_2, \\
-\frac{2b-k_1a}{2b}t^{v_1-k_1-1}+\cdots, & \mbox{if}\ v_1-1=k_2, \\
-t^{k_2-k_1}+\cdots, & \mbox{if}\ v_1-1>k_2,
\end{array} \right. 
$$
with ${\rm I}(g_1,g_2)=k_1$.
\medskip

In the sequel we give an example of a curve $(f)$ for which the analytic type of its polar curve $(af_x+bf_y)$ depends essentially on the direction $(a:b)$.

\begin{ejm} \rm
Consider the curve $(f)$ given parametrically by $(t^5,t^{12}+t^{21})$, that belongs to the eighth family in Example 1. We know that, in this case, $af_x+bf_y=5by^4-10bux^9y-12avx^{11}$, where $u$ and $v$ are units in $\C\{x\}$. This polar is irreducible and is analytically equivalent to a branch with parametrization 
$$\bigg(t^4,t^{11}+t^{14}-\frac 12 t^{17}+\frac{15\sqrt[3]{2}}{2}\big(\frac{12a}{5b}\big)^3t^{21}\bigg).$$
This is a branch of multiplicity four belonging to the fourth Normal Form in Table 3.1. So from the Normal Forms theorem, two such branches corresponding to directions $(a:b)$ and $(a':b')$ with $bb'\neq 0$ are analytically equivalent if and only if one has $\frac{a^3}{b^3}=\frac{a'^3}{b'^3}$. 
\end{ejm}

As a final remark, we refer to \cite{MP} for a rough description of the polars of the members of the equisingularity class determined by the semigroup $\langle 5,11\rangle$, which could be completely described by the methods we exhibited in the present paper.

\end{document}